\documentclass[review]{elsarticle}

\usepackage{lineno,hyperref}
\usepackage{url}
\usepackage{psfrag}
\usepackage{amssymb,amsmath}
\usepackage{amssymb}
\usepackage{dsfont}
\usepackage{multirow}
\usepackage{indentfirst}       
\usepackage{graphicx} 
\usepackage{subfig}           
\usepackage{textcomp}
\usepackage{xspace}
\usepackage{filecontents}
\usepackage{algorithm}
\usepackage{algpseudocode}
\usepackage{color}
\usepackage{tcolorbox}
\usepackage{footnote}
\usepackage{float}
\restylefloat{table}
\usepackage{color}
\usepackage{mathtools}

\newcommand{\removelatexerror}{\let\@latex@error\@gobble}
\newcommand{\RNum}[1]{\uppercase\expandafter{\romannumeral #1\relax}}

\modulolinenumbers[5]

\definecolor{Sahargreen}{rgb}{0.2, 0.5, 0}
\definecolor{Saharblue}{rgb}{0.2, 0.2, 0.8}
\definecolor{Saharred}{rgb}{0.7, 0.05, 0.04}

\journal{EURASIP  signal processing  journal}

\begin{document}
	
	\begin{frontmatter}
		
		\title{Decentralized Decision-Making Over 
			\\ Multi-Task Networks}
		
		\author{Sahar Khawatmi, Abdelhak M. Zoubir\fnref{myfootnote1}}
		\fntext[myfootnote1]{{Member of the European Association for Signal Processing (EURASIP).}}
		\address{Technische Universit\"at
			Darmstadt, Signal Processing Group\\
			64283 Darmstadt, Germany\\
			Email: \{khawatmi, zoubir\}@spg.tu-darmstadt.de}
		\author{Ali H.  Sayed\fnref{myfootnote2}}
		\fntext[myfootnote2]{{Member of the European Association for Signal Processing (EURASIP).}}
		\address{Ecole polytechnique f$\acute{\text{e}}$d$\acute{\text{e}}$rale de Lausanne EPFL,  Adaptive Systems Laboratory\\
		 CH-1015 Lausanne, Switzerland\\
			 Email: ali.sayed@epfl.ch}

	%
		
		
		\begin{abstract}
In important applications involving multi-task networks with multiple objectives, agents in the network need to decide between these multiple objectives and reach an agreement about which single objective to follow for the network. In this work we propose a distributed decision-making algorithm. The agents are assumed to observe data that may be generated by different models. Through localized interactions, the agents reach agreement about which model to track and interact with each other in order to enhance the network performance. We 
investigate the approach for both static and mobile networks. The simulations illustrate the performance of the proposed strategies.
		\end{abstract}
		
		
	\end{frontmatter}
	

\section{Introduction and Related Work}

\noindent  Bio-inspired systems are designed to mimic the 
behavior of  some animal groups such as bee swarms,
birds flying in formation, and schools of 
fish~\cite{bio1,bio4,Avitabile1975,Berg2000,Falko2004,ScotteCamazine2003}.
Diffusion strategies can be used to model some of these coordinated types of behavior, as well as solve inference and   estimation
tasks in a distributed manner over networks~\cite{Ali2013,book2}.
We may distinguish between two types of networks: single-task and multi-task networks. In single-task implementations~\cite{Ali2013,book2}, the networks consist of agents that are interested in the same objective and sense data that are generated by the same model. An analogy would be a school of fish tracking a food source: all elements in the fish school sense distance and direction to the same food source and are interested in approaching it. On the other hand, multi-task networks~\cite{cluster1,cluster2,cluster3,Sahar2,mobile1,new_1,new_2,new_3,talebi2018kalman} involve agents sensing data arising from different models and different clusters of agents may be interested in identifying separate models. A second analogy is a school of fish sensing information about multiple food sources.

In the latter case,     agents 
need to decide between the multiple objectives and reach agreement on following a single objective for the entire network. In the earlier works~\cite{dm,Sahar3},
a  scenario was considered where agents were assumed to sense data arising from {\em two}  models, and a diffusion strategy was developed to enable all agents to agree on estimating a single model. The algorithm developed in~\cite{dm} relied on binary labeling and is applicable only to situations involving two models.  
In this work, we propose an approach for more
than two models.

We  consider a distributed mean-square-error estimation problem over an
$N$-agent network. The connectivity of the agents is described by a graph (see
Fig.~\ref{fig:Partition}).
Data sensed by any particular agent can arise from one of 
different models.   The  objective is to reach an agreement 
among all agents  in the network on one common model to estimate.
Two definitions are introduced:  the
observed model, which refers to the model  from which an agent collects
data,  and the desired model, which refers to the model  the agent decides to 
estimate. The agents do not know which model generated the data they
collect; they also do not know which other agents in
their neighborhood  sense data  arising from the same
model.   Therefore, each agent  needs to determine the subset of
its neighbors that observes the same model. 
This initial step is referred to as {\em clustering}. 
Since the decision-making objective depends on the clustering output, 
errors made during the clustering process have an impact on the global
decision. In this work, we rely on the clustering technique proposed in\cite{Sahar_J} to reduce this effect.

The paper is organized as follows. The network and data model are 
described in Section \RNum{2}. We illustrate the local labeling system and the 
decision-making algorithm in Sections \RNum{3} and \RNum{4}, respectively.  
A special case when the entire  network follows the
model of a specific agent is studied in  Section \RNum{5}. Simulation results 
and  discussion are presented in
Sections \RNum{6} and \RNum{7}, respectively.\\

\noindent {\bf Notation}. We use lowercase letters to denote vectors, uppercase letters for matrices,
plain letters for deterministic variables, and boldface letters for random variables. 
$\mathbb{E}$ denotes the expectation operator and 
$\|\cdot\|$ the Euclidean norm. 
The symbols $\mathds{1}$ and $I$ denote the all-one vector and
identity matrix of appropriate sizes, respectively.  The $k-$th row
(column) of matrix $X$ is denoted by $[X]_{k,:}$ ($[X]_{:,k}$).

\section{Network and Data Model}

\begin{figure}
	\centering\includegraphics[width=20pc]{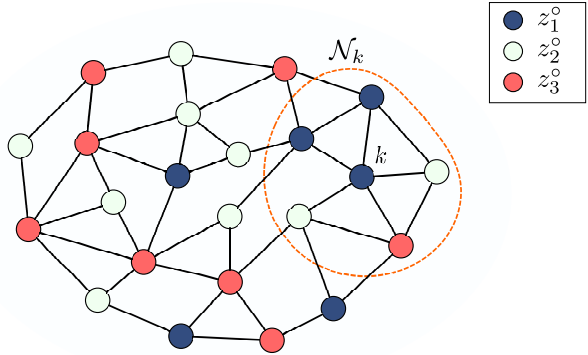}
	\caption{\footnotesize Example of a network topology, agents
		with the same color observe the same model.}
	\label{fig:Partition}
\end{figure}

\noindent Consider a collection of $N$ agents 
distributed  in space, as illustrated in Fig.~\ref{fig:Partition}.
We represent
the network topology  by means of an $N\times N$ adjacency matrix $E$
whose entries $e_{\ell k}$ are defined as follows:   
\begin{equation}
{e}_{\ell k}=\begin{cases} 
1, &  \ell \in \mathcal{N}_{k},\\
0 , &  \textrm{otherwise}
\end{cases}
\end{equation}
where ${\cal N}_{k}$ is  the set of neighbors of agent
$k$ (we denote its size by $n_{k}$). We also write
${\cal N}_k^{-}$ to denote the same neighborhood excluding agent $k$.

Figure~\ref{fig:Partition} shows the network structure where agents
with the same color observe the same model.
We denote the unknown models by $\{z_1^\circ,\ldots,z_C^\circ\}$,
each of size $M\times1$ where $C\leq N$.   Each
agent $k$ observes data generated by one  of  these $C$  unknown models. 
We denote the model observed by agent $k$ by $w_k^\circ$. 
Figure~\ref{fig:Partition} shows that agent $k$  
collects data from model
$z_1^\circ$, in which case  $w_k^\circ= z_1^\circ$. For any other 
agent $\ell$ observing  the same model $z^\circ_1$,  it will hold  that
$w^\circ_{\ell}=z^\circ_1$.
We    stack the $\{w_k^\circ\}$ into a column vector:
\begin{equation}
w^\circ\triangleq \textrm{col}\ \{w_1^\circ,w_2^\circ,\cdots,w_N^\circ\}, \ \ \ w^\circ\in\mathbb{R}^{MN\times 1}.
\end{equation}

At every time instant $i$, every agent $k$ has access to a
scalar measurement $\boldsymbol{d}_k(i)$ and a $1\times M$ regression vector
$\boldsymbol{u}_{k,i}$. The measurements across all  agents are assumed
to be generated via linear regression models of the form:
\begin{equation} \label{eq:data}
\boldsymbol{d}_k(i)=\boldsymbol{u}_{k,i}{w}_{k}^\circ+\boldsymbol{v}_{k}(i).
\end{equation}
All  random processes are assumed to be stationary. Moreover, 
$\boldsymbol{v}_{k}(i)$ is a zero-mean white measurement noise that is  independent over space and has variance
$\sigma_{v,k}^2$. The regression data
$\boldsymbol{u}_{k,i}$ is  assumed to be a zero-mean  random process,
independent over time and space,  and independent of
$\boldsymbol{v}_{\ell}(j)$ for all $k,\ell,i,j$.
We denote the covariance matrix of $\boldsymbol{u}_{k,i}$
by ${R}_{{u},k}=
\mathbb{E} \; \boldsymbol{u}^\intercal_{k,i}\boldsymbol{u}_{k,i}$.  

Agents do not know which model is generating 
their data. 
They also do not know which models are generating 
the data of their neighbors. 
Still, we would like to perform  a learning strategy that  
allows  agents to converge towards one of the models, 
while also learning which of their neighbors share the same model.
Using the
algorithm proposed in~\cite{Sahar_J}, 
each agent $k$ repeats the following steps 
involving an adaptation step 
followed by an aggregation step:
\begin{align}\label{eq:step1}
\boldsymbol{\psi}_{k,i}=& \boldsymbol{\psi}_{k,i-1}+  \mu_k
\boldsymbol{u}^\intercal_{k,i}(  \boldsymbol{d}_k(i) -\boldsymbol{u}_{k,i}
\boldsymbol{\psi}_{k,i-1}) \\ \label{eq:step2}
\boldsymbol{\phi}_{k,i}=& \sum_{\ell=1}^{N}\boldsymbol{a}_{\ell k}(i)
\boldsymbol{\psi}_{\ell,i}
\end{align}
where $\mu_k$ is the
step-size used by agent $k$.  
These updates generate two iterates by 
agent $k$ at time $i$, and which are
denoted by  $\boldsymbol{\psi}_{k,i}$ and $\boldsymbol{\phi}_{k,i}$, 
respectively. We collect the estimated vectors across all agents into the aggregate vector:
\begin{equation}
\boldsymbol{\phi}_i\triangleq \textrm{col}\ 
\{ \boldsymbol{\phi}_{1,i},
\boldsymbol{\phi}_{2,i},\cdots,
\boldsymbol{\phi}_{N,i}\}.
\end{equation}

\noindent  In a manner similar to~\cite{Sahar_J}, we introduce a  clustering
matrix $\boldsymbol{E}_i$. Its structure is similar to the adjacency 
matrix $E$, with
ones and zeros,  except that the value at location
$(\ell,k)$ will be set to one if agent $k$ {\em believes} 
at time instant
$i$ that its neighbor $\ell$ belongs to the same cluster, 
i.e., observes the
same model:
\begin{equation}\label{eq:clsuter_e}
\boldsymbol{e}_{\ell k}(i)=\begin{cases} 
1, &  \textrm{if } \ell\in{\cal N}_k \textrm{ and }  k  \textrm{
	believes that } w_k^\circ=w_\ell^\circ,\\
0 , &  \textrm{otherwise}.
\end{cases}
\end{equation}
These entries help define the neighborhood set
$\boldsymbol{\cal N}_{k,i}$, which consists
of all neighbors at time instant $i$ that  agent $k$
believes share the same model. 
To learn the matrix $\boldsymbol{E}_i$ over time, we apply the clustering technique proposed in~\cite{Sahar_J} to create the estimated clustering matrix $\boldsymbol{F}_i$ of size
$N\times N$ as follows: 
we initialize   $\boldsymbol{\psi}_{k,-1}=0$ and $\boldsymbol{B}_{-1}=\boldsymbol{F}_{-1}=\boldsymbol{E}_{-1}=I_N$.
Where  the matrix $\boldsymbol{B}_{i}$  is of size
$N\times N$.
Each entry  $\boldsymbol{e}_{\ell k}(i)$   is designed using the following
steps from~\cite{Sahar_J}, where $\ell \in {\cal N}_k$:
\begin{equation}     
\boldsymbol{b}_{\ell k}(i)=
\left\{
\begin{aligned}
&1, \text{\ if\ }||\boldsymbol{\psi}_{\ell,i}-\boldsymbol{\phi}_{k,i-1}||^2\leq \alpha  \\          
&0, \text{\ otherwise}
\end{aligned}
\right.
\end{equation}
\begin{equation}  
\boldsymbol{f}_{\ell k}(i)=\nu\times \boldsymbol{f}_{\ell k}(i-1)+(1-\nu)\times \boldsymbol{b}_{\ell k}(i)
\end{equation}
\begin{equation}
\boldsymbol{e}_{\ell k}(i)=\lfloor \boldsymbol{f}_{\ell k}(i)\rceil
\end{equation}
where $\alpha>0$, $0\leq\nu\leq1$,  and the notation $\lfloor
\cdot \rceil$  denotes rounding to the nearest integer.  Using the evolving neighborhoods $\boldsymbol{\cal N}_{k,i}$, the 
entries  $\{\boldsymbol{a}_{\ell k}(i)\}$ in~(\ref{eq:step2})  
are non-negative scalars that
satisfy
\begin{equation}\label{eq:step22}
\boldsymbol{a}_{\ell k}(i) = 0 \ \ \textrm{for}  \ \ell \notin
{\boldsymbol{\mathcal{N}}}_{k,i} , \ \ \ \ \sum_{\ell=1}^N  \boldsymbol{a}_{\ell
	k}(i)=1.
\end{equation}
Although there is a multitude 
of models generating the data that is feeding into
the agents, namely, $\{z_1^\circ, z_2^\circ, \ldots, z_{C}^\circ\}$, 
the objective is to develop a strategy that will allow all agents
to converge towards one of these models.
We refer to this particular choice as the desired model and denote 
it by $z_d^\circ$.

In this way, an agent whose source (observed)  model agrees with 
the desired model, i.e., $w_k^\circ=z^\circ_d$,
will end up tracking its own source. On the other hand, an agent
whose source  model is not
the desired model, i.e.,
$w_k^\circ\neq z^\circ_d$,  will track $z_d^\circ$ instead although it is sensing data generated by a different model.

We define the estimate vector of agent's $k$
desired model by $\boldsymbol{w}_{k,i}$.   
The reason behind  indicating $\boldsymbol{w}_{k,i}$ as
the estimate vector of   \emph{agent's $k$} 
desired model 
instead of the \emph{network's} desired model
is that the agents may have different
desired models before  convergence (steady-state).
Once the  agents reach  agreement among themselves on a single  model, we can then  refer to  $\boldsymbol{w}_{k,i}$ 
as the estimate
vector by agent $k$ of the \emph{network's} desired model.
For the initialization at time
instant $i=1$, each agent assigns
$\boldsymbol{w}_{k,0}=\boldsymbol{\psi}_{k,1}$ (i.e.,
at time instant $i=1$, the
desired model of each agent is a rough estimate of its  own source model).
The decision-making process drives the desired models of all agents to
converge.
For example, if the agents observe 
$C=5$ different
models, the number of the desired models
in the network will decrease with iterations gradually
form  $5$ models down to one model. This is achieved by switching  the
estimate $\boldsymbol{w}_{k,i}$ of some agents 
during the decision-making process 
according to some conditions that are explained later.
However,  agents  do not know which models are desired by  their
neighbors at each time instant $i$.  
Thus, we need to develop a learning
strategy that allows the agents to distinguish the 
individual desired models of their 
neighbors. 

It turns out that in order for the objective of the network  
to be met, 
it is important for agents to combine the estimates of 
their neighbors in a judicious manner because, 
unbeknown to 
the agents,   some of their neighbors may be wishing to estimate 
different models. 
If cooperation is performed blindly with all neighbors, 
then performance can deteriorate  with agents 
converging  to non-existing locations.
For this reason, and motivated by the discussion
from~\cite{dm}, we add 
the step~(\ref{eq:w}) below  after~(\ref{eq:step1}) and~(\ref{eq:step2}), 
which involves two sets of
combination coefficients from two matrices 
$\boldsymbol{\dot{A}}_i$ and
$\boldsymbol{\ddot{A}}_i$. There are two main   ideas behind the
construction~(\ref{eq:w}). First, it is meant to  let  each agent $k$  cooperate
only with the subset of  neighbors that
share the same desired model as it does.
Second, it also lets  each agent $k$ combine
$\boldsymbol{\phi}_{\ell,i}$ if the desired model of 
agent $k$ at time instant
$i$ is the same as $\ell$'s  observed model:
\begin{equation}\label{eq:w}
\boldsymbol{w}_{k,i}=\sum_{\ell =1}^N \boldsymbol{\dot{a}}_{\ell
	k}(i) \boldsymbol{\phi}_{\ell,i}
+\sum_{\ell =1}^N \boldsymbol{\ddot{a}}_{\ell k}(i)
\boldsymbol{w}_{\ell,i-1}.
\end{equation}
Note that the matrices $\boldsymbol{\dot{A}}_i$ and
$\boldsymbol{\ddot{A}}_i$ are not constructed from matrix
$\boldsymbol{{A}}_i$. 
The selection of the  non-negative   
coefficients $\{\boldsymbol{\dot{a}}_{\ell k}(i)\}$
and  $\{\boldsymbol{\ddot{a}}_{\ell k}(i)\}$ is 
explained in Section \RNum{4}.

We summarize the main five  steps  of the approach:
\begin{enumerate}
	\item  {\em Learning the observed models of the neighbors}. This step  is
	performed by building the matrix $\boldsymbol{E}_i$ 
	in step~(\ref{eq:clsuter_e}).  
	The  information provided
	by each entry $\boldsymbol{e}_{\ell k}(i)$  
	is  whether the corresponding agents $\ell$ and $k$ have the
	same observed model or not. 
	\item   {\em Learning and labeling the desired model
		of the neighbors at each time
		instant $i$}.
	This step allows the agents to distinguish the 
	individual desired models of their 
	neighbors at  time
	instant $i$.  The  information provided  by this step is  
	the number of different 
	models that are desired by neighbors and how many times each  model is 
	repeated at  time
	$i$ among neighbors.
	\item   {\em Decision-making step} by switching the desired model 
	of some agents to let
	the network converge to only one model.
	\item  {\em Learning the desired models of the neighbors after the switching step}.
	This step  is performed by building the matrix  $\boldsymbol{H}_i$ 
	in step~(\ref{eq:h}) in Section \RNum{4}.
	The  information provided
	by each entry $\boldsymbol{h}_{\ell k}(i)$  
	is  whether the corresponding agents $\ell$ and $k$ have the
	same desired model or not after the switching step.  
	\item   {\em Updating the estimate vectors} $\{\boldsymbol{w}_{k,i}\}$ 
	by sharing data  thoughtfully with the subset of the neighbors that
	share the same desired model.
\end{enumerate}

 \section{Local Labeling}

\begin{figure}
	\centering\includegraphics[width=18pc]{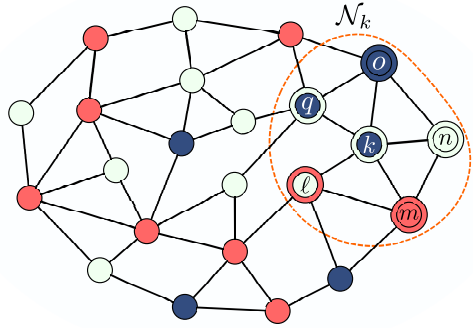}
	\caption{\footnotesize Example of an agent $k$ and its neighborhood ${\cal
			N}_k$. The inner color indicates
		the observing model while the outer one indicates the current desired model.}
	\label{fig_Ch5_2:Labeling}
\end{figure}

\noindent   Each agent needs to learn the desired models of its
neighbors to proceed with the decision-making
process and let the network converge to only one model. In this step, instead of
only estimating  whether two agents have the same desired model or not, 
the construction involves a local
labeling procedure that enables every agent
to estimate in real-time how many different 
models are desired by its neighborhood. 

For this purpose, we associate with 
\emph{each} agent $k$ an $n_k \times n_k$ matrix
$\boldsymbol{Y}_i^k$  with  
entries $\{\boldsymbol{y}^k_{\ell m}(i)\}$   given by:
\begin{equation}\label{eq:y}
\boldsymbol{y}^k_{\ell m}(i)=\begin{cases} 
1, &   \textrm{if }\| \boldsymbol{w}_{m,i-1} - \boldsymbol{w}_{\ell,i-1}\|^2
\leq \beta,\\ 
0 , &  \textrm{otherwise}
\end{cases}
\end{equation} 
for some small threshold $\beta>0$. 
Whenever $\boldsymbol{y}^k_{\ell m}(i)=1$,
agent $k$ believes at time  instant $i$ that  its
neighbors $\ell$ and $m$ wish to estimate the 
same  desired model.  	
On account of that,   the variables
$\boldsymbol{w}_{m,i-1}$ and 
$\boldsymbol{w}_{\ell,i-1}$ which are used in the
test~(\ref{eq:y}) are presenting
the current desired model of  agents $m$
and $\ell$, respectively.
It is clear from~(\ref{eq:y}) that the matrix 
$ \boldsymbol{Y}^k_i$ is  symmetric and has  
ones on the diagonal. An example is depicted in Fig.~\ref{fig_Ch5_2:Labeling} 
where agents having the same inner
color observe the  same model, while the outer color
indicates the  model in which the agent is interested
(or towards which the
agent is moving in mobile networks).
The corresponding matrix $\boldsymbol{Y}^k_i$ has the following entries:
\begin{equation}\label{eq:Y}
\boldsymbol{Y}^k_i=
{  \begin{array}{c}
	k \\
	\ell \\
	m\\
	n\\
	o\\
	q 
	\end{array}}
\overbracket[0.001pc][0.001pc]{\left[
	\begin{array}{cccccc} 
	\color{Sahargreen} 1 & \color{Saharred} 0 & \color{Saharred}0 & \color{Sahargreen}1 & \color{Saharblue}0
	& \color{Sahargreen} 1 \\
	\color{Sahargreen}  0 & \color{Saharred}1 & \color{Saharred}1 & \color{Sahargreen} 0 & \color{Saharblue}0
	& \color{Sahargreen} 0 \\
	\color{Sahargreen}  0 & \color{Saharred}1 & \color{Saharred}1 &\color{Sahargreen}  0 & \color{Saharblue}0
	& \color{Sahargreen} 0 \\
	\color{Sahargreen}  1 & \color{Saharred}0 & \color{Saharred}0 & \color{Sahargreen} 1 & \color{Saharblue}0
	& \color{Sahargreen} 1 \\
	\color{Sahargreen}  0 & \color{Saharred}0 & \color{Saharred}0 & \color{Sahargreen} 0 & \color{Saharblue}1
	& \color{Sahargreen} 0 \\
	\color{Sahargreen} 1 & \color{Saharred}0 & \color{Saharred}0 & \color{Sahargreen} 1 & \color{Saharblue}0
	& \color{Sahargreen} 1 \\
	\end{array} \right]}^{  \ \ k  \ \  \ \ell  \ \  \ m \  \ \  n  \ \ 
	\ o \ \ \ q \ \  }.
\end{equation} 
From~(\ref{eq:Y}) agents that share  the same desired model will  have
identical columns in matrix $\boldsymbol{Y}^k_i$, namely, if agents
$m$ and $\ell$ have the same desired model at time instant $i$, this implies
that:
$[\boldsymbol{Y}^k_i]_{:,m}=[\boldsymbol{Y}^k_i]_{:,\ell}$.
We denote the local label of each agent $\ell \in {\cal
	N}_k$ by agent $k$ as $\boldsymbol{l}^k_{\ell}(i)$. The local  label $\boldsymbol{l}^k_{\ell}(i)$
is updated at each time instant $i$ using the following relation:
\begin{equation}\label{eq:label}
\boldsymbol{l}^k_{\ell}(i)={ \cal B}([\boldsymbol{Y}^k_i]_{:,\ell})
\end{equation}
where ${\cal B}(\cdot)$ is  a function that converts the input sequence from
binary to decimal.  For the   example in~(\ref{eq:Y}), we have
\begin{align*}
&\boldsymbol{l}^k_{k}(i) ={\cal B}(100101)=37,\\
&\boldsymbol{l}^k_{\ell}(i) ={\cal B}(011000)=24,\\
&\boldsymbol{l}^k_{m}(i) ={\cal B}(011000)=24,\\
& \boldsymbol{l}^k_{n}(i) ={\cal B}(100101)=37,\\
& \boldsymbol{l}^k_{o}(i) ={\cal B}(000010)=2,\\
& \boldsymbol{l}^k_{q}(i) ={\cal B}(100101)=37.
\end{align*}
We define  the number of  
desired models within ${\cal N}_k$  at time instant $i$ by
$\boldsymbol{C}_{k}(i)$.
After updating matrix $\boldsymbol{Y}^k_i$  and generating the local labels
$\{\boldsymbol{l}^k_{\ell}(i)\}$, agent $k$ counts how many models are desired
by its neighborhood to  update
$\boldsymbol{C}_{k}(i)$. 
In the example~(\ref{eq:Y}), 
agent $k$  distinguishes at time instant $i$  three  desired
models $\{2,24,37\}$, i.e., $\boldsymbol{C}_{k}(i)=3$.  Agent $k$  labels 
these three different
models locally as:
$\{2,24,37\}$.

In addition, agent $k$ determines which model among these 
$\boldsymbol{C}_{k}(i)$
models has the maximum number of   followers. A   \emph{follower} of
a model is an agent  that wishes to estimate and track this model.  
We define the largest set of agents 
belonging to ${\cal N}_k$ and
following  the same desired model at time  instant $i$  by
$\boldsymbol{ \cal Q}_{k,i}$.
In the example, agent $k$ assigns the majority set
at time instant $i$ as follows:   $\boldsymbol{\cal
	Q}_{k,i}=\{k,n,q\}$ which has the label $37$ and 
is repeated three times among
other labels.

\section{Decision-Making Over Multi-Task Networks}

\begin{figure}
	\centering\includegraphics[width=18pc]{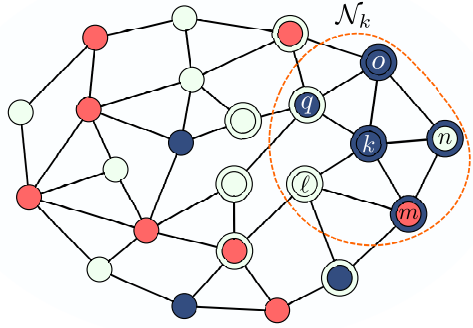}
	\caption{\footnotesize Example of the equilibrium case. All agents within ${\cal 
			N}_k$ belong  to the majority sets among their neighbors.}
	\label{fig:equil}
\end{figure}

\noindent Using the information provided by matrix $\boldsymbol{Y}^k_i$, agent
$k$ can capture how many agents within its neighbors
follow the same desired
model at time instant $i$. Once agent  $k$ and all its 
neighbors agree on a single desired model, the  matrix
$\boldsymbol{Y}^k_i$ will end up being of the following form with unit entries everywhere:
\begin{equation}\label{eq:agree}
\boldsymbol{Y}^k_i=
{  \begin{array}{c}
	k \\
	\ell \\
	m\\
	n\\
	o\\
	q 
	\end{array}}
\overbracket[0.001pc][0.001pc]{\left[
	\begin{array}{cccccc} 
	1 & 1& 1& 1 & 1
	&  1 \\
	1 & 1& 1& 1 & 1
	&  1 \\
	1 & 1& 1& 1 & 1
	&  1 \\
	1 & 1& 1& 1 & 1
	&  1 \\
	1 & 1& 1& 1 & 1
	&  1 \\
	1 & 1& 1& 1 & 1
	&  1 \\
	\end{array} \right]}^{  \ \ k  \ \  \ \ell  \ \  \ m \  \ \  n  \ \ 
	\ o \ \ \ q \ \  }.
\end{equation} 
We define the degree of  agreement by each agent $k$ among
its neighbors ${\cal N}_k$ as
\begin{equation}\label{eq:p}
\boldsymbol{p}_{k}(i)=\frac{[\boldsymbol{Y}^k_i]_{k,:}
	\mathds{1}}{n_k}.
\end{equation}
Equally, having $\boldsymbol{p}_{k}(i)=1$  means that agent $k$ and all  of its neighbors  have agreed on a common  desired model. On the other hand, if 
$\boldsymbol{p}_{k}(i)\neq1$, then the following switching step is applied:
\begin{equation}\label{eq:switch}
\boldsymbol{w}_{k,i-1}\leftarrow
\begin{cases}
\boldsymbol{w}_{\ell,i-1}, &  
\textrm{if } k\notin  \boldsymbol{ \cal Q}_{k,i}  \textrm{ for any }  \; \ell  \in
\boldsymbol{ \cal Q}_{k,i},\\
\boldsymbol{w}_{n,i-1} , &    \textrm{if } k\in  \boldsymbol{ \cal Q}_{k,i}
\textrm{ and }  \boldsymbol{C}_{k}(i)=2,  \ \forall n \in { \cal N}_{k},
\\
\boldsymbol{w}_{k,i-1} , &  \textrm{otherwise}.
\end{cases}
\end{equation}
The main idea of the  switching step is for each agent $k$ to make a new
decision or to keep the previous one. The first case of~(\ref{eq:switch})
implies that agent $k$ does not belong to the majority desired model set
$\boldsymbol{ \cal Q}_{k,i}$ at time instant $i$. Therefore, agent $k$ changes its decision and switches into the desired model of the majority set $\boldsymbol{ \cal Q}_{k,i}$. The second
case in~(\ref{eq:switch}) is applied to prevent an unwanted equilibrium situation. This problem may arise 
when only two desired models remain in
${\cal N}_k$.  In this case, if  all agents in ${\cal N}_k$
belong  to the majority set,  this leads to a situation in which
no agent  in ${\cal N}_k$ will change its decision anymore.
An example is shown in Fig.~\ref{fig:equil}  where the outer color of the
agents indicate the desired model. We indicate  only the desired
model of agent's $k$
neighbors and their neighbors.
Figure~\ref{fig:equil} shows that all
agents within ${\cal  N}_k$ belong  to a  majority set and no agent 
in ${\cal N}_k$ will
change its decision anymore, e.g. agents $q$ and $\ell$ 
belong to the  majority set
among their neighbors, as well as agents $k$, $m$, $n$, and $o$. Namely,
\begin{align*}
&k\in  \boldsymbol{ \cal Q}_{k,i},   
m\in  \boldsymbol{ \cal Q}_{m,i},   
n\in  \boldsymbol{ \cal Q}_{n,i}, \textrm{and }
o\in  \boldsymbol{ \cal Q}_{o,i}  \ \  (\textrm{with } z_1^\circ),\\
&\ell\in  \boldsymbol{ \cal Q}_{\ell,i} \textrm{ and } q\in  \boldsymbol{ \cal
	Q}_{q,i}  \ \ (\textrm{with } z_2^\circ).
\end{align*}
To break the equilibrium,   an agent  that recognizes
these two models  picks randomly one of
the two desired models.

From~(\ref{eq:switch}), we can conjecture that 
the network will probably
converge to the most observable model, since the initial desired model by each agent is
its own observed model. This fact remains true even with 
the random
switching in the second case of~(\ref{eq:switch}), 
because in that case  the more
repeated desired model within ${\cal N}_k$ has
the highest probability to be picked.

To proceed with the
cooperation and sharing information among the agents within the subset that has
the same desired model at time instant $i$, we define an
$N\times N$ matrix $\boldsymbol{H}_i$.
The coefficients $\{\boldsymbol{h}_{\ell k}(i)\}$ are
updated after the switching step~(\ref{eq:switch}) using 
a   test that is quite similar to~(\ref{eq:y})
and is  applied between each agent $k$ and its neighbors as follows:
\begin{equation}\label{eq:h}
\boldsymbol{h}_{\ell k}(i)=\begin{cases}
1, &   \textrm{if }\| \boldsymbol{w}_{k,i-1} - \boldsymbol{w}_{\ell,i-1}\|^2
\leq \beta,\\
0 , &  \textrm{otherwise}.
\end{cases}
\end{equation}
According to  matrix $\boldsymbol{H}_i$,
each agent knows which subset of
its neighbors has the same desired model as it does   after the
switching step at time instant $i$.
Having $\boldsymbol{h}_{\ell k}(i)=1$  means that
$\ell$ and $k$ have the same desired model at time instant $i$.
We define an $N\times N$ combination matrix $\boldsymbol{G}_i$ as follows:
\begin{equation}\label{eq:g}
\boldsymbol{G}_i={\cal F}(\boldsymbol{H}_i)
\end{equation}
where  ${\cal F}(\cdot)$ is some
function which satisfies
\begin{equation}\label{eq:step33}
\boldsymbol{g}_{\ell k}(i) = 0  \ \ \textrm{ if } \ \ \boldsymbol{h}_{\ell
	k}(i)=0, \ \ \ \ \sum_{\ell=1}^N  \boldsymbol{g}_{\ell k}(i)=1
\end{equation}
An example of ${\cal F}(\cdot)$ is the uniform function which generates a left-stochastic matrix $\boldsymbol{G}_i$	where each entry  $\boldsymbol{g}_{\ell k}(i)$ is given by
\begin{equation}\label{eq_Ch2:A_Uniform}
\boldsymbol{g}_{\ell k}(i)=\begin{cases} 
\frac{1}{ \sum_{n=1}^{N}\boldsymbol{h}_{n
		k}(i)}, &  \textrm{if } \boldsymbol{h}_{\ell
	k}(i)=0, \\
0 , &  \textrm{otherwise}.
\end{cases}
\end{equation}
Matrix $\boldsymbol{G}_i$ by itself does  not have enough  information  for
proceeding and updating the estimate $\boldsymbol{w}_{k,i}$. The agents
still need  knowledge about which data to be combined from each neighbor.
Therefore, matrix $\boldsymbol{G}_i$ is split into two matrices
$\boldsymbol{\dot{A}}_i$  and $\boldsymbol{\ddot{A}}_i$.
The weight of the entry $ \boldsymbol{g}_{\ell k}(i)$ goes to
$\boldsymbol{\dot{a}}_{\ell k}(i)$ if the 
desired model of agent $k$ at time instant $i$ is the same as $\ell$'s
observed model.
Otherwise, $\boldsymbol{\ddot{a}}_{\ell k}(i)$   obtains the weight $
\boldsymbol{g}_{\ell k}(i)$. 
The coefficients $\{\boldsymbol{\dot{a}}_{\ell
	k}(i)\}$ and $\{\boldsymbol{\ddot{a}}_{\ell k}(i)\}$ 
for $\ell \in {\cal N}_k$ are
updated using the following steps:
\begin{align}\label{eq:dot_a}
\boldsymbol{\dot{a}}_{\ell k}(i)&=\begin{cases}
\boldsymbol{g}_{\ell k}(i), &   \textrm{if }\| \boldsymbol{w}_{k,i-1} -
\boldsymbol{\psi}_{\ell,i}\|^2 \leq \beta,\\
0 , &  \textrm{otherwise}.
\end{cases} \\ \label{eq:ddot_a}
\boldsymbol{\ddot{a}}_{\ell k}(i)&=\begin{cases}
\boldsymbol{g}_{\ell k}(i), &    \textrm{if } \boldsymbol{\dot{a}}_{\ell
	k}(i)=0, \\
0 , &  \textrm{otherwise}.
\end{cases}
\end{align}
In~(\ref{eq:dot_a}),  the case that  $\boldsymbol{\psi}_{\ell,i}$ is close to
$\boldsymbol{w}_{\ell,i-1}$ implies that the observed model of agent $\ell$ is
the same as the  desired model of agent $k$ at time instant $i$.
The estimate
$\boldsymbol{w}_{k,i}$ is updated using~(\ref{eq:w}).
Algorithm 1 summarizes  the decision-making scheme.
\begin{algorithm}[]
	\begin{algorithmic}
		\State Initialize
		$\boldsymbol{A}_{0}=\boldsymbol{\dot{A}}_{0}=\boldsymbol{\ddot{A}}_{0}=\boldsymbol{E}_{0}=\boldsymbol{H}_{0}=\boldsymbol{G}_{0}=I$
		\State Initialize $\boldsymbol{\psi}_{0}=\boldsymbol{\phi}_{0}=0$ and
		$\boldsymbol{p}_0=0$
		\For {$i> 0$}
		\vspace{0.2cm}
		\For {$k=1,\ldots,N$}
		\vspace{-0.1cm}
		\begin{align}
		\boldsymbol{\psi}_{k,i} = \boldsymbol{\psi}_{k,i-1}+  \mu_k
		\boldsymbol{u}^\intercal_{k,i}(  \boldsymbol{d}_k(i) -\boldsymbol{u}_{k,i}
		\boldsymbol{\psi}_{k,i-1})
		\end{align}
		\State assign $\boldsymbol{w}_{k,0}=\boldsymbol{\psi}_{k,1}$ at $i=1$
		\State update $\{\boldsymbol{a}_{\ell k}(i)\}$ according
		to~(\ref{eq:step22})   
		\begin{align}
		\boldsymbol{\phi}_{k,i}= \sum_{\ell=1}^{N}\boldsymbol{a}_{\ell k}(i)
		\boldsymbol{\psi}_{\ell,i}
		\end{align}
		\State  update $\boldsymbol{Y}^k_i$ according to~(\ref{eq:y}) 
		\State find $\boldsymbol{\cal Q}_{k,i}$ and $\boldsymbol{C}_k(i)$ 
		\State update $\boldsymbol{p}_{k}(i)$ according to~(\ref{eq:p})
		\If{$\boldsymbol{p}_{k}(i)\neq1$}
		\State switch $ \boldsymbol{w}_{k,i-1}$ according  to~(\ref{eq:switch}) 
		\State resend $ \boldsymbol{w}_{k,i-1}$
		\EndIf
		\vspace{0.2cm}
		\For {$\ell \in {\cal N}_{k}$}
		\vspace{0.2cm}
		\State update $\{\boldsymbol{h}_{\ell k}(i)\}$ according
		to~(\ref{eq:h})
		\State    update $\{\boldsymbol{g}_{\ell k}(i)\}$ according
		to~(\ref{eq:g})
		\State update $\{\boldsymbol{\dot{a}}_{\ell k}(i)\}$ according
		to~(\ref{eq:dot_a})
		\State update $\{\boldsymbol{\ddot{a}}_{\ell k}(i)\}$  according
		to~(\ref{eq:ddot_a})
		\vspace{0.2cm}
		\EndFor
		\vspace{-0.07cm}
		\begin{equation}
		\boldsymbol{w}_{k,i}=\sum_{\ell =1}^N \boldsymbol{\dot{a}}_{\ell
			k}(i) \boldsymbol{\phi}_{\ell,i}
		+\sum_{\ell =1}^N \boldsymbol{\ddot{a}}_{\ell k}(i)
		\boldsymbol{w}_{\ell,i-1}
		\end{equation}
		\EndFor 
		\vspace{0.05cm}
		\EndFor 
	\end{algorithmic}
	\caption{(Decentralized decision-making  scheme)}
	\label{alg:S}
\end{algorithm}

\section{Following the Observed Model of a Specific Agent}

\begin{figure}
	\centering\includegraphics[width=18pc]{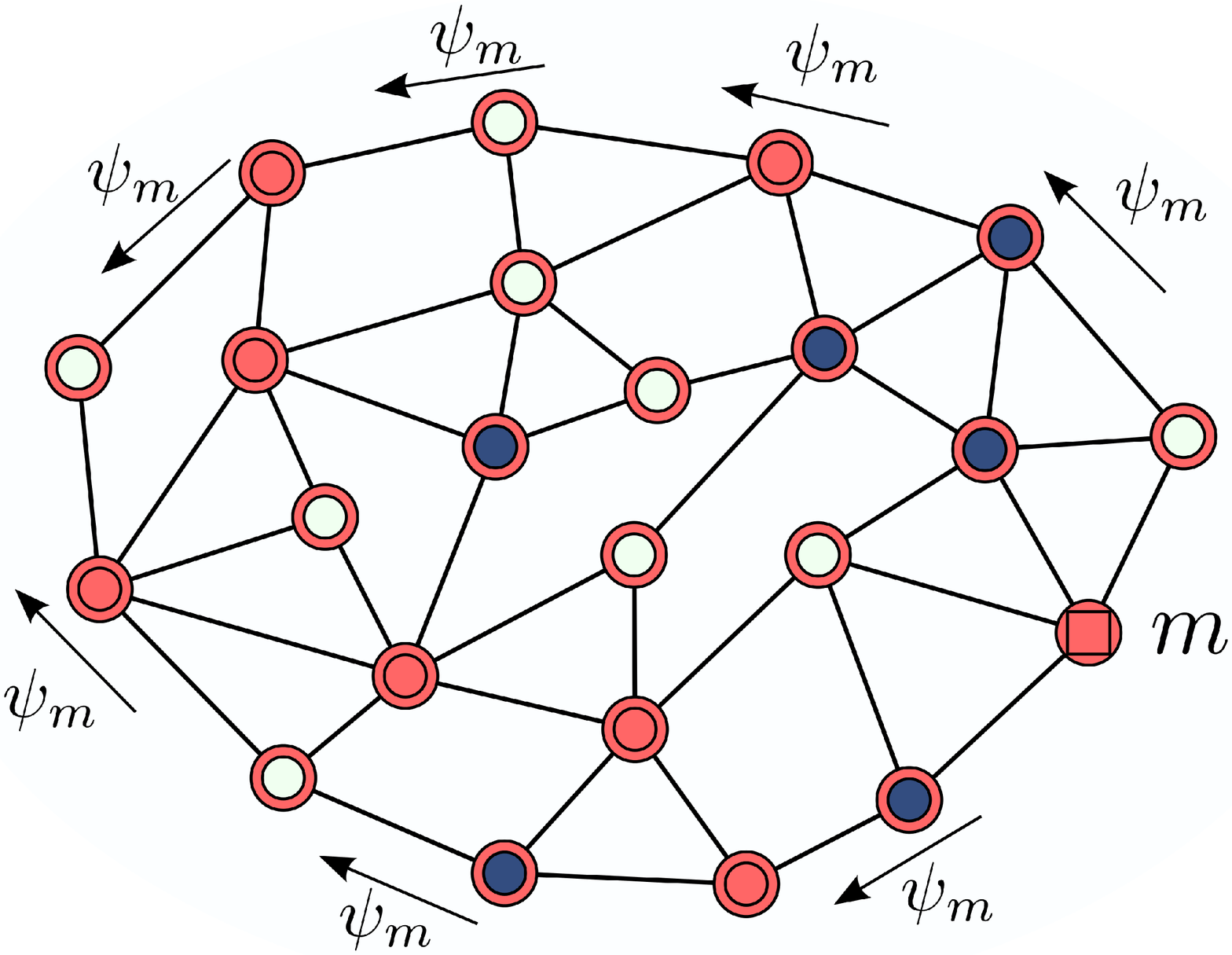}
	\caption{\footnotesize Final  decision of a network after following the
		model of the specific agent $m$. The inner color indicates
		the observing model while the outer one indicates the desired model.
		The arrows represent the spreading process of $\boldsymbol{\psi}_{m,i}$  through the
		network.}
	\label{fig:Agent}
\end{figure}

\noindent   In this section the goal is to let the whole network follow the
observed model of  some specific agent  $m$, as  shown in
Fig.~\ref{fig:Agent} where agent $m$ observes  model $z_3^\circ$ (red),
therefore, the network 
converges in a distributed manner to estimate the
model $z_3^\circ$. The first
step is to spread the $\boldsymbol{\psi}_{m,i}$
among agents and keep updating
it over time. This step aims at having a  copy (reference) of
$\boldsymbol{\psi}_{m,i}$ by all agents in the network.
Agents  keep updating the copy of
$\boldsymbol{\psi}_{m,i}$ for two reasons. First, to have a more accurate version of the vector 
$\boldsymbol{\psi}_{m,i}$,
which indicates the desired  model of the network. Second,  to endow the algorithm to work in  non-stationary situations, if 
drift is happening in agent $m$'s model.

\begin{figure*}
	\centering\includegraphics[width=28pc]{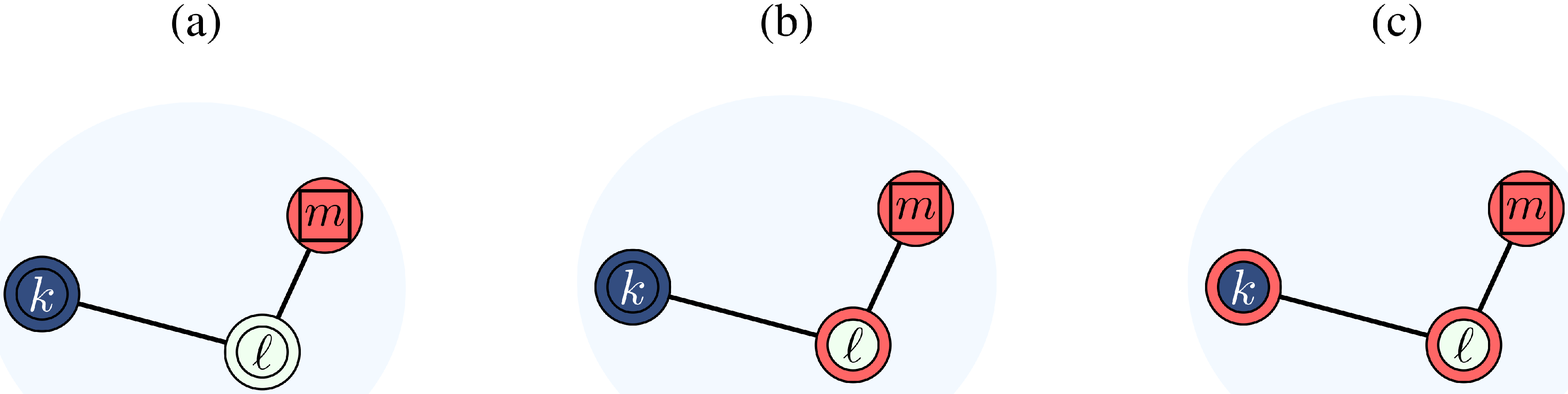}
	\vspace{0.5cm}
	\caption{\footnotesize   Example of the spreading process of
		$\boldsymbol{\psi}_{m,i}$  from agent $m$ to agent $k$ over time. The inner
		color indicates the observing model while the outer one indicates the desired model.}
	\label{fig:Agent2}
\end{figure*}

We denote the  copy vector of $\boldsymbol{\psi}_{m,i}$  by
agent $k$ 
by $\boldsymbol{\breve{\psi}}_{k,i}$ and refer to it as  the \emph{anchor
	vector}. 
Agents are informed beforehand about the index $m$ of the specific agent  that
they should follow. If $m\in {\cal N}_k$, this implies that agent $k$  receives
the anchor vector directly from agent $m$. If not, i.e., $m\notin {\cal N}_k$,
then agent $k$ depends on  another agent $\ell\in {\cal
	N}_k$ that has already a copy of $\boldsymbol{\psi}_{m,i}$. 
Agent  $k$ stores the index of this source agent. 
The index of the source agent of agent $k$
is denoted by $\boldsymbol{s}_k(i)$.
Note that the anchor vector
$\boldsymbol{\breve{\psi}}_{k,i}$ is not the final 
estimate of the desired model. 

The circulation process of $\boldsymbol{\psi}_{m,i}$ in a
distributed manner needs  cooperation among agents.
In case that agent $k$ has no direct
link to receive   data from agent $m$, i.e., $m \not\in {\cal N}_{k}$,
agent $k$
gets one of the $\boldsymbol{\breve{\psi}}_{\ell,i-1}$ provided that
$\boldsymbol{s}_\ell(i)\neq0$. If $\boldsymbol{s}_\ell(i)\neq0$
this implies that agent $\ell$ has
already a source to
update its $\boldsymbol{\breve{\psi}}_{\ell,i}$,
regardless whether $m \in {\cal
	N}_{\ell}$ or not. In other words,  $\boldsymbol{s}_\ell(i)\neq0$
means that agent $\ell$
finds a direct or indirect link to agent $m$. Therefore,
it is important for each agent $k$ to store the agent's index of  its source.  An example is 
shown in Fig.~\ref{fig:Agent2} where $m \in {\cal N}_\ell$ but 
$m \notin {\cal
	N}_k$. First,  the anchor vectors and the source
agents for agents $k$ and $\ell$ at time instant $i=0$
(Fig.~\ref{fig:Agent2}(a)) are given, respectively, by
\begin{align}
\boldsymbol{\breve{\psi}}_{k,0}=0,  \boldsymbol{s}_k(0)=0,
\boldsymbol{\breve{\psi}}_{\ell,0}=0,  \boldsymbol{s}_\ell(0)=0.
\end{align}
The anchor vectors and the source
agents for agents $k$ and $\ell$ at time instants $i=\{1,2\}$
(Fig.~\ref{fig:Agent2}(b) and (c)) are given, respectively, by
\begin{align}
&\boldsymbol{\breve{\psi}}_{k,1}=0,  \boldsymbol{s}_k(1)=0,
\boldsymbol{\breve{\psi}}_{\ell,1}=  \boldsymbol{{\psi}}_{m,1}
,   \boldsymbol{s}_\ell(1)=m,\\
&\boldsymbol{\breve{\psi}}_{k,2}=\boldsymbol{\breve{\psi}}_{\ell,1}, 
\boldsymbol{s}_k(2)=\ell,
\boldsymbol{\breve{\psi}}_{\ell,2}= \boldsymbol{{\psi}}_{m,2}
,  \boldsymbol{s}_\ell(2)=m.
\end{align}
Agents update  their anchor vectors 
$\{\boldsymbol{\breve{\psi}}_{k,i}\}$
at each time instant $i$ by the following step:
\begin{equation}\label{eq:switch2}
\boldsymbol{\breve{\psi}}_{k,i}=
\begin{cases}
\boldsymbol{{\psi}}_{m,i}, &   \textrm{if } m \in {\cal N}_{k},\\
\boldsymbol{\breve{\psi}}_{\ell,i-1}, &  \textrm{if } \ell \in {\cal N}_{k}
\textrm{ and }  \boldsymbol{s}_k(i)=0
\textrm{ and }  \boldsymbol{s}_\ell(i)\neq0,\\
\boldsymbol{\breve{\psi}}_{\ell,i-1}, &  \textrm{if } \ell \in {\cal N}_{k}
\textrm{ and }  \boldsymbol{s}_k(i)=\ell,  \\
\boldsymbol{\breve{\psi}}_{k,i-1}, &   \textrm{otherwise}
\end{cases}
\end{equation}
where $\boldsymbol{\breve{\psi}}_{m,i}=\boldsymbol{{\psi}}_{m,i}$  for
agent $m$ itself. The source of the 
anchor vector is updated simultaneously as follows:
\begin{equation}\label{eq:source}
\boldsymbol{s}_k(i)=
\begin{cases}
m, &   \textrm{if } m \in {\cal N}_{k},\\
\ell, &  \textrm{if }   \boldsymbol{s}_k(i)=0
\textrm{ and }  \boldsymbol{s}_\ell(i)\neq0,\\
\boldsymbol{s}_k(i-1), &   \textrm{otherwise}.
\end{cases}
\end{equation}
Similarly to the previous section, the next
step  is to update the coefficients $\{\boldsymbol{h}_{\ell k}(i)\}$ using
the following test:
\begin{equation}\label{eq:h2}
\boldsymbol{h}_{\ell k}(i)=\begin{cases}
1, &   \textrm{if } \boldsymbol{s}_\ell(i)\neq0  \textrm{ and
} \boldsymbol{s}_k(i)\neq0,\\
0 , &  \textrm{otherwise}.
\end{cases}
\end{equation}
Again, having $\boldsymbol{s}_k(i)\neq0$
leads to the situation that agent $k$ has the anchor vector and has been
informed about the decision of the network, therefore,  agent $k$ can start sharing
information with the other agents whose $\boldsymbol{s}_\ell(i)\neq0$ as well to
estimate the desired model.
The matrix $\boldsymbol{G}_i$ will be generated using~(\ref{eq:g}).
Agents  update the  coefficients of  both matrices
$\boldsymbol{\dot{A}}_i$ and $\boldsymbol{\ddot{A}}_i$  using the following steps:
\begin{align}\label{eq:dot_a2}
\boldsymbol{\dot{a}}_{\ell k}(i)&=\begin{cases}
\boldsymbol{g}_{\ell k}(i), &   \textrm{if }\| \boldsymbol{\breve{\psi}}_{k,i}
- \boldsymbol{\psi}_{\ell,i}\|^2 \leq \beta,\\
0 , &  \textrm{otherwise}.
\end{cases} \\ \label{eq:ddot_a2}
\boldsymbol{\ddot{a}}_{\ell k}(i)&=\begin{cases}
\boldsymbol{g}_{\ell k}(i), &    \textrm{if } \boldsymbol{\dot{a}}_{\ell
	k}(i)=0, \\
0 , &  \textrm{otherwise}.
\end{cases}
\end{align}
Then, the estimate $\boldsymbol{w}_{k,i}$ is  updated using Eq.~(\ref{eq:w}).
According to~(\ref{eq:dot_a2}) and~(\ref{eq:w}), agent $k$
combines $\boldsymbol{\phi}_{\ell,i}$ if the desired
model of the network (which is represented by the anchor vector
$\boldsymbol{\breve{\psi}}_{k,i}$ of agent $k$) is close to the observed model
of agent $\ell$ that is represented by
$\boldsymbol{{\psi}}_{\ell,i}$.
Algorithm 2 summarizes  the steps of the  approach for
following the observed model
of a specific agent $m$.
\begin{algorithm}[]
	\begin{algorithmic}
		\State Initialize
		$\boldsymbol{A}_{0}=\boldsymbol{\dot{A}}_{0}=\boldsymbol{\ddot{A}}_{0}=\boldsymbol{E}_{0}=\boldsymbol{H}_{0}=\boldsymbol{G}_{0}=I$
		\State Initialize
		$\boldsymbol{\psi}_{0}=\boldsymbol{\breve{\psi}}_{0}=\boldsymbol{\phi}_{0}=0$
		and $\boldsymbol{s}_0=0$
		\vspace{0.2cm}
		\For {$i> 0$}
		\vspace{0.1cm}
		\For {$k=1,\ldots,N$}
		\vspace{-0.1cm}
		\begin{align}
		\boldsymbol{\psi}_{k,i} = \boldsymbol{\psi}_{k,i-1}+  \mu_k
		\boldsymbol{u}^\intercal_{k,i}(  \boldsymbol{d}_k(i) -\boldsymbol{u}_{k,i}
		\boldsymbol{\psi}_{k,i-1})
		\end{align}
		\State assign $\boldsymbol{w}_{k,0}=\boldsymbol{\psi}_{k,1}$ at $i=1$
		\State update $\{\boldsymbol{a}_{\ell k}(i)\}$ according
		to~(\ref{eq:step22})  
		\begin{align}
		\boldsymbol{\phi}_{k,i}= \sum_{\ell=1}^{N}\boldsymbol{a}_{\ell k}(i)
		\boldsymbol{\psi}_{\ell,i}
		\end{align}
		\State update $\boldsymbol{\breve{\psi}}_{k,i}$ according to~(\ref{eq:switch2})
		\State update $\boldsymbol{s}_k(i)$ according to~(\ref{eq:source})
		\vspace{0.2cm}
		\For {$\ell \in {\cal N}_{k}$}
		\vspace{0.2cm}
		\State update $\{\boldsymbol{h}_{\ell k}(i)\}$ according
		to~(\ref{eq:h2})
		\State    update $\{\boldsymbol{g}_{\ell k}(i)\}$ according
		to~(\ref{eq:g})
		\State update $\{\boldsymbol{\dot{a}}_{\ell k}(i)\}$ according
		to~(\ref{eq:dot_a2})
		\State update $\{\boldsymbol{\ddot{a}}_{\ell k}(i)\}$  according
		to~(\ref{eq:ddot_a2})
		\vspace{0.2cm}
		\EndFor
		\vspace{-0.07cm}
		\begin{equation}
		\boldsymbol{w}_{k,i}=\sum_{\ell =1}^N \boldsymbol{\dot{a}}_{\ell
			k}(i) \boldsymbol{\phi}_{\ell,i}
		+\sum_{\ell =1}^N \boldsymbol{\ddot{a}}_{\ell k}(i)
		\boldsymbol{w}_{\ell,i-1}
		\end{equation}
		\EndFor 
		\vspace{0.1cm}
		\EndFor 
	\end{algorithmic}
	\caption{(Following the observed model of a specific agent)}
	\label{alg:S}
\end{algorithm}

\section{Simulation Results and Discussion}

\subsection{Static Network}

\noindent  We consider a connected network with 80 randomly distributed
agents. The agents observe data originating from $C=3$ different models.
Each model $z_j^\circ\in \mathbb{R}^{M\times1}$ is generated as
follows: $z_j^\circ=[ {r_1},\ldots,{r_M}]^\intercal$ where
${r_m}\in[1,-1]$, $M=2$. The assignment of the agents to models is
random.  The maximum
number of neighbors is $n_{k}=7$.  We set $\{\alpha,\beta,\nu,\mu\}
=\{0.04,0.08,0.005,0.01\}$. We use the uniform combination policy to generate
the coefficients $\{\boldsymbol{a}_{\ell k}(i)\}$ and $\{\boldsymbol{g}_{\ell
	k}(i)\}$.

\begin{figure}
	\psfrag{V}[c]{\small $\sigma_{v,k}^2$}
	\centerline{\includegraphics[width=9cm]{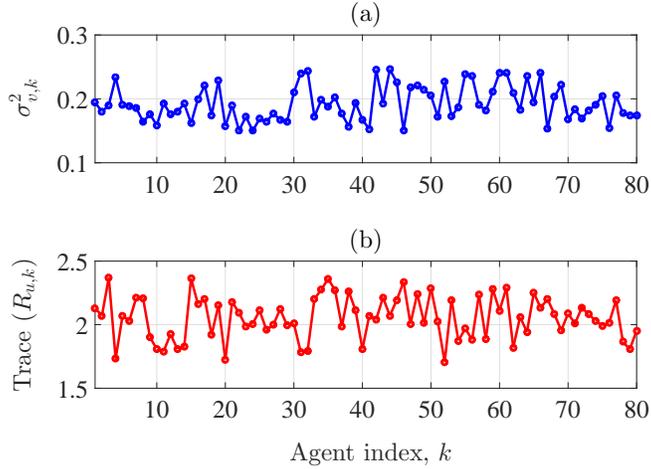}}
	\caption{ \footnotesize Statistical noise and signal profiles over the
		network.}
	\label{fig:N}
\end{figure}

\begin{figure}
	\vspace{-0.0cm}
	\centering
	\centerline{\includegraphics[width=24pc]{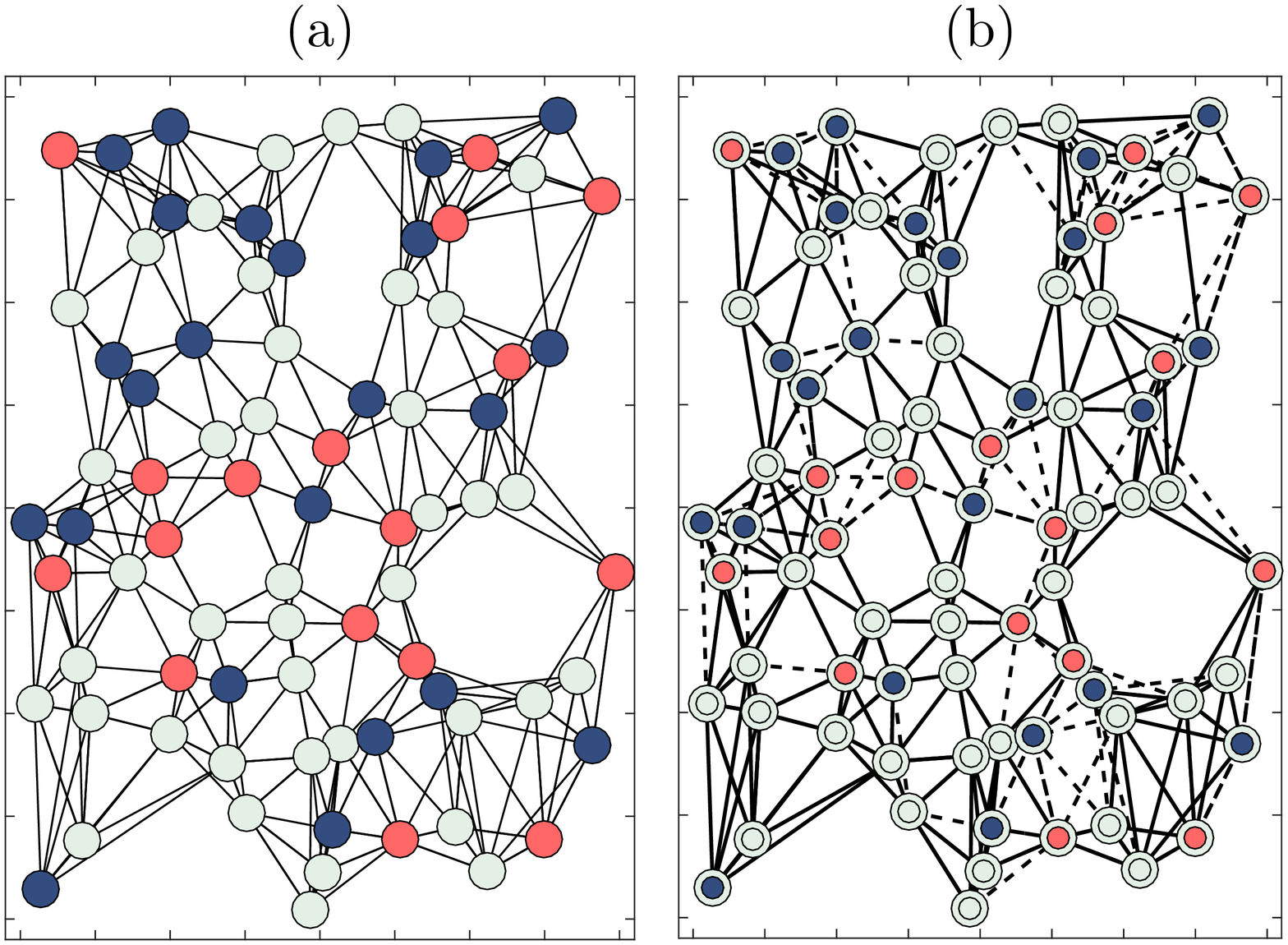}}
	\caption{ \footnotesize Network topology (a) and
		final decision of the agents where the bold (dashed) links
		represent $\{\boldsymbol{\dot{a}}(i)\}$ ($\{\boldsymbol{\ddot{a}}(i)\}$)
		at steady-state (b).}
	\label{fig:Topology}
\end{figure}

Figure~\ref{fig:N} shows the statistical profile of the
regressors and noise across the agents. The regressors are of size
$M=2$  zero-mean Gaussian, independent in time and space, and have
diagonal covariance matrices ${R}_{u,k}$.
Figure~\ref{fig:Topology} shows the topology of one    of 100 Monte
Carlo experiments.
Agents having the same inner color observe the  same model, while the outer
color indicates the desired model at steady-state.
\begin{figure}
	\centering 
	\centerline{\includegraphics[width=22pc]{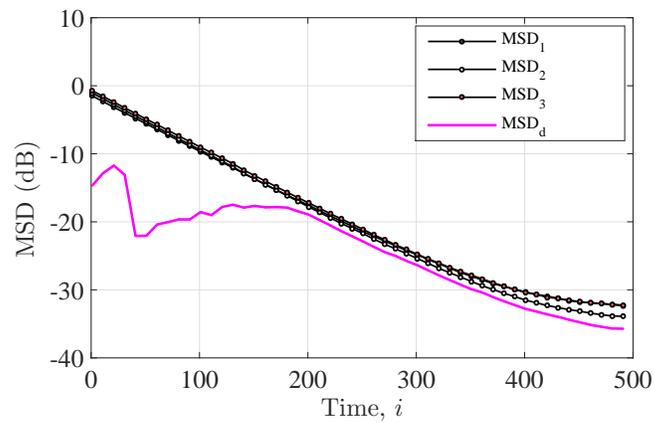}}
	\caption{ \footnotesize Transient mean-square deviation (MSD).}
	\label{fig:D_No}
\end{figure}

The transient network mean-square deviation
(MSD) regarding
each observed model $z^\circ_j$ at each time instant $i$ is defined by
\begin{equation}
\textrm{MSD}_{j}(i)\triangleq\frac{1}{|{\cal C}_j|}\sum_{k \in {\cal C}_j}
||z^\circ_j-{\boldsymbol{\phi}_{k,i}}||^2
\end{equation}
where $j=1,\ldots,C$ and each $ \textrm{MSD}_{j}$ is computed for  agents
belonging to ${\cal C}_j$. The transient network mean-square deviation (MSD) for
the whole network regarding the desired model at each time instant $i$ is
defined by
\begin{equation}
\textrm{MSD}_\textrm{d}(i)\triangleq\frac{1}{N}\sum_{k=1}^{N}
||z^\circ_d-{\boldsymbol{w}_{k,i}}||^2
\end{equation}
where $z_d^\circ$ is the desired model when the whole network
agrees on one common desired model, i.e., $\textrm{MSD}_\textrm{d}(i)$ is
only computed at the instants when all $\{\boldsymbol{p}_k(i)\}=1$.
Figure~\ref{fig:D_No} depicts the simulated transient mean-square
deviation (MSD) of the network for all observed models and for the network
desired model.
%
Table~\ref{table1} displays the success rate of
the decision-making to agree on one model for different numbers of
observed models, $C\in\{2,3,4,5\}$. The proposed strategy appears to provide good  success rate under the simulated conditions. 
\begin{table}[H]
	\centering
	\begin{tabular}{c|c|c|c|c|c}
		$C$  & 2 & 3 & 4 & 5  \\
		\hline
		\small Success rate   & 99\% & 98\%  &  99\% &  99\%
	\end{tabular}
	\caption {\small Decision-making success rate
		for different  $C$.}
	\label{table1}
\end{table}

\begin{figure}
	\vspace{-0.1cm}
	\centering
	\centerline{\includegraphics[width=22pc]{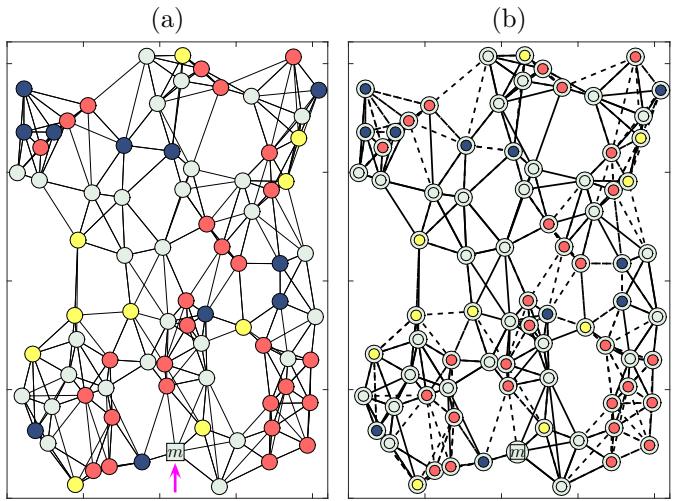}}
	\caption{ \footnotesize Network topology (a) and
		final decision of the agents to follow the model of agent $m$ where the
		bold (dashed) links represent $\{\boldsymbol{\dot{a}}(i)\}$ ($\{\boldsymbol{\ddot{a}}(i)\}$)
		at steady-state (b).}
	\label{fig:Topology_Nodek}
\end{figure}
\begin{figure}
	\vspace{-0cm}
	\centering 
	\centerline{\includegraphics[width=22pc]{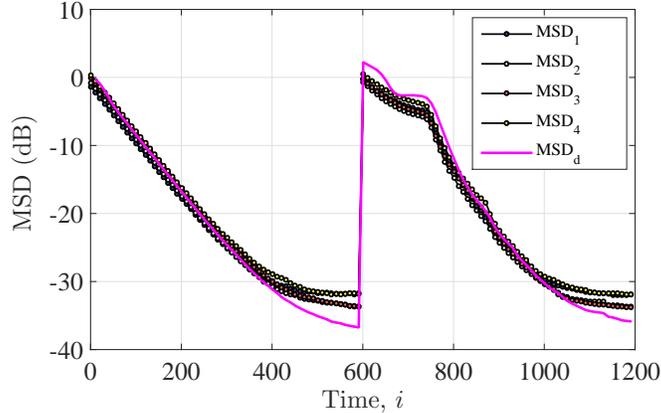}}
	\caption{ \footnotesize Transient mean-square deviation (MSD).}
	\label{fig:MSD}
\end{figure}

Regarding the application of
following the observed model of a specific agent $m$,
Fig.~\ref{fig:Topology_Nodek} shows the topology of one case from  100
different experiments.
Agents are observing $C=4$ different models.
Agent $m=10$, which is  represented  by a square, is the specific
agent whose observed model  the whole network wishes to follow.
Figure~\ref{fig:MSD}  shows the transient mean-square deviation MSD of 100
different experiments when a change in the model assignment occurs
suddenly at time instant $i=600$. The success rate of
the decision-making to agree on the observed model of
agent $m$ was observed to be 100\% in this simulation.

\subsection{Mobile Network}
\begin{figure}
	\vspace{-0.3cm}
	\psfrag{V}[cB]{\small $\sigma_{v,k}^2$}
	\centerline{\includegraphics[width=9cm]{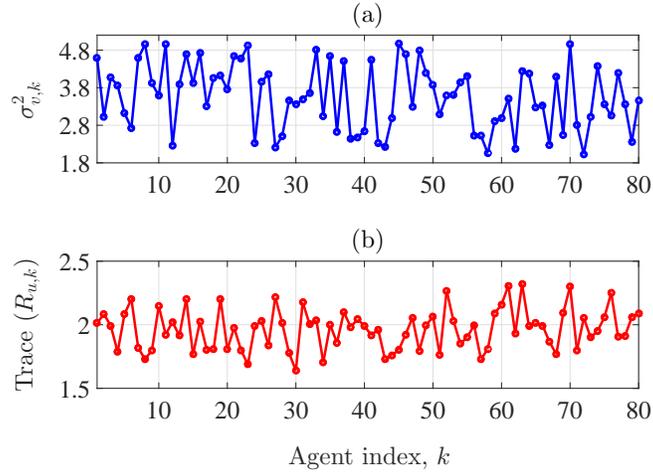}}
	\caption{ \footnotesize Statistical noise and signal profiles over the
		mobile network.}
	\label{fig:N_M}
\end{figure}

\noindent  We consider a network with 80 randomly
distributed mobile agents~\cite{bio4}.
The agents observe data originating
from four different models (sources) $C=4$, where
$w_{r_m}\in[50,-50]$. The objective of the network is to have all agents
track and move towards  only one model (source).
Figure~\ref{fig:N_M} shows the
statistical profile of the regressors and noise across the agents. Every agent
$k$ updates its location  according
to the motion mechanism described  in~\cite{Sahar3}.

\begin{figure}[]
	\vspace{-0.1cm}
	\hspace{-0.2cm}
	\centerline{\includegraphics[width=23pc]{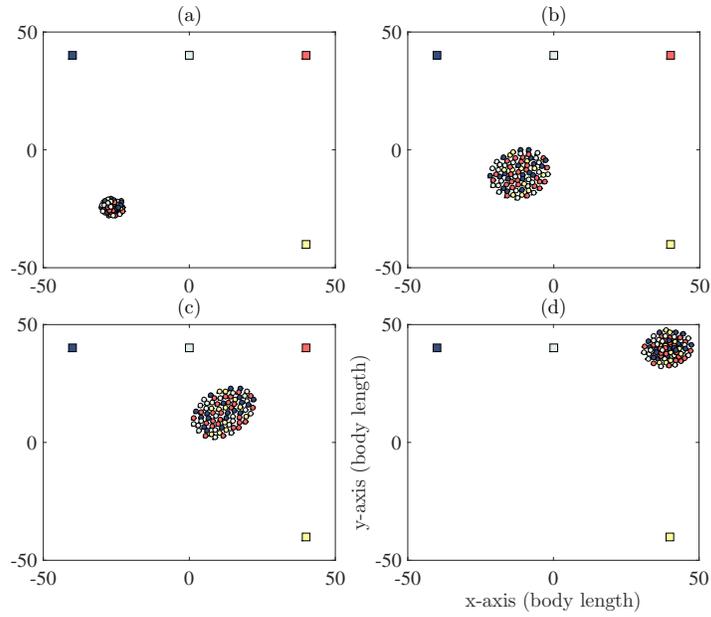}}
	\caption{ \footnotesize Maneuver of the agents with four  sources over
		time (a) $i$=1, (b) $i$=200, (c) $i$=500, and (d) $i$=1000. The
		unit  length  is the body length of a agent.}
	\label{fig:Maneuver}
\end{figure}

\begin{figure}
	\vspace{-0.2cm}
	\centering 
	\centerline{\includegraphics[width=22pc]{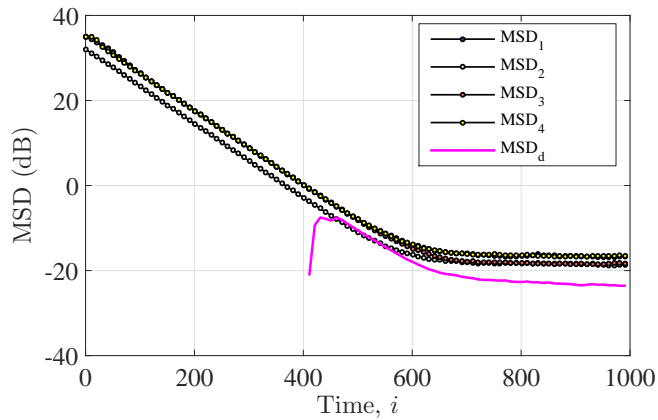}}
	\caption{\footnotesize Transient mean-square deviation (MSD) of the mobile
		network.}
	\label{fig:Mobile}
\end{figure}

Figure~\ref{fig:Maneuver}
shows the maneuver of the agents over time where the models (sources) are
represented by  squares.
Figure~\ref{fig:Mobile} represents the  transient network mean-square
deviation (MSD)  obtained
by averaging over 100 independent Monte Carlo experiments.

\section{Conclusion}

\noindent We  have proposed a distributed algorithm 
that allows   agents over 
multi-task networks to follow only one common model while
proceeding with the
estimation process. Agents use a local labeling step to  distinguish 
the multiple desired models of their  neighbors. 
Simulation results illustrate the operation of the algorithms and its performance. 

\section*{Acknowledgements}

\noindent The research for this paper was financially supported by Technische Universit\"at
Darmstadt, Signal Processing Group and NSF grant CCF-1524250.

\end{document}